arXiv:1407.2818v1 [math.AP] 10 Jul 2014

# ON THE LOW MACH NUMBER LIMIT OF COMPRESSIBLE FLOWS IN EXTERIOR MOVING DOMAINS

EDUARD FEIREISL, ONDŘEJ KREML, VÁCLAV MÁCHA AND ŠÁRKA NEČASOVÁ

Abstract. We study the incompressible limit of solutions to the compressible barotropic Navier-Stokes system in the exterior of a bounded domain undergoing a simple translation. The problem is reformulated using a change of coordinates to fixed exterior domain. Using the spectral analysis of the wave propagator, the dispersion of acoustic waves is proved by the means of the RAGE theorem. The solution to the incompressible Navier-Stokes equations is identified as a limit.

Institute of Mathematics of the Academy of Sciences of the Czech Republic
Žitná 25, 115 67 Praha 1, Czech Republic

## 1. Introduction

We study the compressible barotropic Navier-Stokes equations

$$\partial_t \varrho + \operatorname{div}_x(\varrho \mathbf{u}) = 0 \tag{1.1}$$

$$\partial_t(\varrho \mathbf{u}) + \operatorname{div}_x(\varrho \mathbf{u} \otimes \mathbf{u}) + \nabla_x p(\varrho) = \operatorname{div}_x \mathbb{S}(\nabla_x \mathbf{u}), \tag{1.2}$$

with $\varrho(t,x)$ denoting the density of the fluid and $\mathbf{u}(t,x)$ denoting the velocity of the fluid being the unknowns of the system. The pressure $p(\varrho)$ is a given function and $\mathbb{S}(\nabla_x \mathbf{u})$ denoting the viscous stress tensor is given by

$$\mathbb{S}(\nabla_x \mathbf{u}) = \mu \left( \nabla_x \mathbf{u} + (\nabla_x \mathbf{u})^T - \frac{2}{3} \operatorname{div}_x \mathbf{u} \mathbb{I} \right) + \eta \operatorname{div}_x \mathbf{u} \mathbb{I}, \qquad \mu > 0, \quad \eta \geq 0. \tag{1.3}$$

We consider the system (1.1)-(1.2) in the exterior of a bounded domain undergoing a translation. Therefore we introduce the following notation. Let $T > 0$ and let $m(t) : [0,T] \mapsto \mathbb{R}^3$ be a given smooth function satisfying $m(0) = 0$. Let moreover $\Sigma_0 \subset \mathbb{R}^3$ be a fixed bounded domain of the class $C^2$. We define $\Sigma_t := \{x \in \mathbb{R}^3, x - m(t) \in \Sigma_0\}$. Let us denote moreover $\Omega_t := \mathbb{R}^3 \setminus \Sigma_t$ and the space-time cylinder $Q_t := \{(s,x) \in [0,t] \times \mathbb{R}^3, x \in \Omega_s\}$. The system (1.1)-(1.2) has to be satisfied in $Q := Q_T$.

We denote by $\Gamma_t$ the boundary of the domain $\Omega_t$ and denote $\Gamma := \{(t,x) \in (0,T) \times \mathbb{R}^3, x \in \Gamma_t\}$. Moreover for $(t,x) \in \Gamma$ let $\nu(t,x)$ be a unit outer normal to $Q$ in space-time. Finally, we write $\nu(t,x) = (\nu_t(t,x), \mathbf{n}(t,x))$. We complement the system (1.1)-(1.2) with boundary conditions expressing impermeability of the boundary and complete slip of the fluid on the boundary

$$\mathbf{u} \cdot \mathbf{n} = -\nu_t, \quad [\mathbb{S}(\nabla_x \mathbf{u}) \cdot \mathbf{n}] \times \mathbf{n} = 0, \quad \text{on } \Gamma. \tag{1.4}$$

All authors acknowledge the support of the GAČR (Czech Science Foundation) project GA13-00522S in the general framework of RVO: 67985840.

Moreover we prescribe the behavior of the fluid at infinity

$$\varrho \to \overline{\varrho} > 0, \quad \mathbf{u} \to 0 \quad \text{as } |x| \to \infty. \tag{1.5}$$

The pressure $p \in C[0,\infty) \cap C^1(0,\infty)$ is assumed to fulfill

$$p(0) = 0, \quad p'(\varrho) > 0 \text{ for all } \varrho > 0, \quad \lim_{\varrho\to\infty} \frac{p'(\varrho)}{\varrho^{\gamma-1}} = p_\infty > 0 \text{ for a certain } \gamma > \frac{3}{2}. \tag{1.6}$$

In the regime, where the speed of sound dominates the characteristic speed of the fluid, the system (1.1)-(1.2) can be rescaled using dimensionless quantities. Assuming the Mach number to be of order $\varepsilon$, whereas all other dimension numbers to be of order 1, we obtain the scaled Navier-Stokes system

$$\partial_t \varrho + \mathrm{div}_x\, \varrho \mathbf{u} = 0 \tag{1.7}$$

$$\partial_t(\varrho \mathbf{u}) + \mathrm{div}_x(\varrho \mathbf{u} \otimes \mathbf{u}) + \frac{1}{\varepsilon^2} \nabla_x p(\varrho) = \mathrm{div}_x\, \mathbb{S}(\nabla_x \mathbf{u}). \tag{1.8}$$

We consider a family of weak solutions $(\varrho_\varepsilon, \mathbf{u}_\varepsilon)$ to the system (1.7)-(1.8) with boundary conditions (1.4)-(1.5) emanating from the initial data

$$\varrho_\varepsilon(0,\cdot) = \varrho_{0,\varepsilon} = \overline{\varrho} + \varepsilon \varrho^{(1)}_{0,\varepsilon}, \quad \mathbf{u}_\varepsilon(0,\cdot) = \mathbf{u}_{0,\varepsilon}, \tag{1.9}$$

where

$$\|\varrho^{(1)}_{0,\varepsilon}\|_{L^2(\Omega_0)} + \|\varrho^{(1)}_{0,\varepsilon}\|_{L^\infty(\Omega_0)} \le c, \ \mathbf{u}_{0,\varepsilon} \to \mathbf{u}_0 \text{ weakly in } L^2(\Omega_0;\mathbb{R}^3). \tag{1.10}$$

Our aim is to prove that in a certain sense the weak solutions of the system (1.7)-(1.8) converge to weak solutions of the incompressible Navier-Stokes system:

$$\overline{\varrho}\partial_t \mathbf{U} + \overline{\varrho}\, \mathrm{div}_x(\mathbf{U} \otimes \mathbf{U}) + \nabla_x \Pi = \mathrm{div}_x\, \mathbb{S}(\nabla_x \mathbf{U}), \tag{1.11}$$

$$\mathrm{div}_x\, \mathbf{U} = 0, \tag{1.12}$$

together with boundary conditions (1.4) and initial condition

$$\mathbf{U}(0,\cdot) = \mathbf{H}(\mathbf{u}_0) \quad \text{in } \Omega_0,$$

where $\mathbf{H}$ denotes the Helmholtz projection to the space of solenoidal functions in $\Omega_0$.

Although the problem is very simple in the sense that the underlying rigid object undergoes only translations and the fluid slips on its boundary, the authors are not aware of any rigorous mathematical result on the low Mach number limit for the *ill prepared* initial data (1.9), (1.10). Originated by the pioneering result of Lions [14] on the *existence* of large data weak solutions for the compressible Navier-Stokes system, Desjardins et al. [4], [5] Lions and Masmoudi [15] (see also the surveys Danchin [3], Masmoudi [16], Schochet [18] and the references cited therein) employed the framework of weak solutions to singular incompressible limits for problems confined to fixed spatial domains. Similar problems on a bounded time dependent domain, with *prescribed* boundary motion, have been studied only recently in [8].

Similarly to [7], in order to show compactness of the convective term in the momentum equation, we use the dispersive estimates for the underlying acoustic equation. To this end, the problem is transformed to a fixed spatial domain giving rise to a perturbed acoustic equation, for which the desired decay estimates follow from the energy bounds combined with an application of the celebrated RAGE theorem.

The paper is organized as follows. In Section 2, we collect some preliminary material concerning the weak solutions of both primitive and target system and state our main result. Section 3 contains the basic estimates derived directly from the associated energy balance. As a consequence, we deduce weak convergence of the solutions of the scaled system in Section 4. The acoustic equation is derived in Section 5, and the dispersive estimates obtained in Section 6. The proof of the main result is then completed in Section 7.

## 2. Preliminaries and main result

2.1. **Weak solutions of the primitive system.** We start with a simple Lemma.

**Lemma 2.1.** *There exists a function* $\mathbf{V} : Q \mapsto \mathbb{R}^3$ *such that* $\mathbf{V} \in C^1([0,T], W^{1,2} \cap W^{1,\infty}(\Omega_t))$ *satisfying*

$$\begin{aligned} \operatorname{div}_x \mathbf{V} &= 0 \quad \textit{in } Q \\ \mathbf{V} &= -\nu_t \mathbf{n} \quad \textit{on } \Gamma. \end{aligned} \tag{2.1}$$

*and* $\mathbf{V}(t,x) \equiv 0$ *for all* $t \in [0,T]$ *and* $|x| > R$ *for sufficiently large* $R$.

*Proof.* For any $R > 0$ denote $B_R^c := \{x, |x| > R\}$. Take $R > 0$ such that the set $B_R^c \subset \Omega_t$ for all $t \in [0,T]$. Observing that for all $t \in [0,T]$ it holds

$$m'(t) \cdot \mathbf{n} = \nu_t \text{ on } \Gamma_t, \ \int_{\Gamma_t} m'(t) \cdot \mathbf{n} \, \mathrm{d}S = 0,$$

we can define $\mathbf{V}(t,x)$ as the Bogovskii solution (see Bogovskii [1], Galdi [11, Chapter 3]) to the problem

$$\begin{aligned} \operatorname{div}_x \mathbf{V}(t,\cdot) &= 0 \quad \text{in } \Omega_t \cap B_R \\ \mathbf{V}(t,\cdot) &= -\nu_t \mathbf{n} \quad \text{on } \Gamma_t \\ \mathbf{V}(t,\cdot) &= 0 \quad \text{on } B_R^c. \end{aligned} \tag{2.2}$$

□

Note that in particular $(\mathbf{u} - \mathbf{V}) \cdot \mathbf{n} = 0$ on $\Gamma$. Now we are ready to define weak solutions to the primitive system.

**Definition 2.2.** We say that a couple $(\varrho, \mathbf{u})$ is a weak solution to the compressible Navier-Stokes system (1.7)-(1.8) with boundary conditions (1.4)-(1.5) and initial conditions (1.9) if

1. $\varrho - \overline{\varrho} \in L^\infty(0,T,(L^2 + L^\gamma)(\Omega_t))$, $\mathbf{u} \in L^2(0,T,W^{1,2}(\Omega_t))$
2. The continuity equation (1.7) is satisfied in a weak sense and in a renormalized form,i.e.
$$\begin{aligned} &\int_{\Omega_t} (\varrho(t) + b(\varrho(t)))\varphi(t,\cdot) \, \mathrm{d}x - \int_{\Omega_0} (\varrho_{0,\varepsilon} + b(\varrho_{0,\varepsilon}))\varphi(0,\cdot) \, \mathrm{d}x \\ &= \int_{Q_t} \big( (\varrho + b(\varrho))\partial_t \varphi + (\varrho + b(\varrho))\mathbf{u} \cdot \nabla_x \varphi + (b(\varrho) - b'(\varrho)\varrho) \operatorname{div}_x \mathbf{u} \varphi \big) \, \mathrm{d}x \, \mathrm{d}t \end{aligned} \tag{2.3}$$
holds for all $t \in [0,T]$, all $\varphi \in C_c^\infty(\overline{Q})$, such that $\varphi(T,\cdot) = 0$ in $\Omega_T$, and any $b \in C^1[0,\infty)$ such that $b(0) = 0$ and $b'(r) = 0$ for large $r$.
3. Balance of linear momentum is satisfied in a weak sense, i.e.
$$\begin{aligned} &\int_{\Omega_t} (\varrho \mathbf{u})(t) \cdot \varphi(t) \, \mathrm{d}x - \int_{\Omega_0} (\varrho_{0,\varepsilon} \mathbf{u}_{0,\varepsilon}) \cdot \varphi(0) \, \mathrm{d}x \\ &= \int_{Q_t} \left( \varrho \mathbf{u} \cdot \partial_t \varphi + \varrho(\mathbf{u} \otimes \mathbf{u}) : \nabla_x \varphi + \frac{1}{\varepsilon^2} p(\varrho) \operatorname{div}_x \varphi - \mathbb{S}(\nabla_x \mathbf{u}) : \nabla_x \varphi \right) \mathrm{d}x \, \mathrm{d}t \end{aligned} \tag{2.4}$$

holds for all $t \in [0,T]$ and all $\varphi \in C_{cn} := \{\varphi \in C_c^\infty(\overline{Q}), \varphi \cdot \mathbf{n} = 0 \text{ on } \Gamma\}$.

(4) Energy inequality holds, i.e.

$$\int_{\Omega_t} \left[\frac{1}{2}\varrho\,|\mathbf{u}|^2 + \frac{1}{\varepsilon^2}\Big(P(\varrho) - P'(\overline{\varrho})(\varrho - \overline{\varrho}) - P(\overline{\varrho})\Big)\right](t,\cdot)\,\mathrm{d}x + \int_{Q_t} \mathbb{S}(\nabla_x \mathbf{u}) : \nabla_x \mathbf{u}\,\mathrm{d}x\,\mathrm{d}t$$
$$\leq \int_{\Omega_0} \left[\frac{1}{2}\rho_{0,\varepsilon}\,|\mathbf{u}_{0,\varepsilon}|^2 + \frac{1}{\varepsilon^2}\Big(P(\varrho_{0,\varepsilon}) - P'(\overline{\varrho})(\varrho_{0,\varepsilon} - \overline{\varrho}) - P(\overline{\varrho})\Big)\right]\mathrm{d}x + \int_{\Omega_t} (\varrho\mathbf{u}\cdot\mathbf{V})(t,\cdot)\,\mathrm{d}x$$
$$- \int_{\Omega_0} \rho_{0,\varepsilon}\mathbf{u}_{0,\varepsilon}\cdot\mathbf{V}(0,\cdot)\,\mathrm{d}x + \int_{Q_t} \left(\mathbb{S}(\nabla_x\mathbf{u}) : \nabla_x\mathbf{V} - \varrho\mathbf{u}\otimes\mathbf{u} : \nabla_x\mathbf{V} - \varrho\mathbf{u}\cdot\partial_t\mathbf{V}\right)\,\mathrm{d}x\,\mathrm{d}t, \tag{2.5}$$

where

$$P(\varrho) = \varrho\int_1^\varrho \frac{p(z)}{z^2}\,\mathrm{d}z.$$

The following existence result of weak solutions to the compressible Navier-Stokes system in moving domains was proved in [8].

**Theorem 2.3.** *Let $\varepsilon > 0$. There exists a weak solution $(\varrho_\varepsilon, \mathbf{u}_\varepsilon)$ to the system* (1.7)-(1.8) *with boundary conditions* (1.4)-(1.5) *and initial conditions* (1.9) *in the sense of Definition 2.2.*

Strictly speaking, the result of [8] covers the case of a *bounded* physical space, however, the extension to the exterior problem is straightforward, see Sýkora [19].

Now we define weak solutions to the target system, the incompressible Navier–Stokes equations.

**Definition 2.4.** We say that $\mathbf{U}$ is a weak solution to (1.11)-(1.12) in $Q$ if

(1) $\mathbf{U} \in L^\infty(0,T,L^2(\Omega_t)) \cap L^2(0,T,W^{1,2}(\Omega_t))$
(2) $\operatorname{div}_x \mathbf{U} = 0$ a.a. in $(0,T)\times\Omega_t$
(3) $\mathbf{U}\cdot\mathbf{n} = \mathbf{V}\cdot\mathbf{n}$ on $\Gamma$
(4) For all $\varphi \in C_c^\infty([0,T)\times\Omega_t)$, $\operatorname{div}_x\varphi = 0$, $\varphi\cdot\mathbf{n} = 0$ on $\Gamma$, it holds

$$\int_Q \left(\overline{\varrho}\mathbf{U}\cdot\partial_t\varphi + \overline{\varrho}(\mathbf{U}\otimes\mathbf{U}) : \nabla_x\varphi\right)\,\mathrm{d}x\,\mathrm{d}t$$
$$= \int_Q \mathbb{S}(\nabla_x\mathbf{U}) : \nabla_x\varphi\,\mathrm{d}x\,\mathrm{d}t + \int_{\Omega_0}\overline{\varrho}\mathbf{U}_0\cdot\varphi(0,\cdot)\,\mathrm{d}x \tag{2.6}$$

We are now ready to state the main theorem.

**Theorem 2.5.** *Let $m(t) \in C^2([0,T])$. Let $(\varrho_\varepsilon, \mathbf{u}_\varepsilon)$ be a sequence of weak solutions to the compressible Navier-Stokes system* (1.7)–(1.8) *with boundary conditions* (1.4)–(1.5) *and initial conditions* (1.9)*. Then, at least for a suitable subsequence*

$$\operatorname{ess\,sup}_{t\in(0,T)}\|\varrho_\varepsilon(t,\cdot) - \overline{\varrho}\|_{(L^2+L^\gamma)(\Omega_t)} \leq C\varepsilon \tag{2.7}$$

$$\mathbf{u}_\varepsilon \;\to\; \mathbf{U} \text{ in } L^2(K)$$

*for any compact $K \subset Q$, where $\mathbf{U}$ is a weak solution to the incompressible Navier-Stokes system* (1.11)–(1.12) *with initial data $\mathbf{U}_0 = \mathbf{H}[\mathbf{u}_0]$.*

The rest of the paper is devoted to the proof of Theorem 2.5.

## 3. Uniform estimates

Following [10] , we introduce the *essential part*

$$[f_\varepsilon]_{ess} := f_\varepsilon 1_{\{\frac{\overline{\varrho}}{2}<\varrho_\varepsilon<2\overline{\varrho}\}},$$

and the *residual part*

$$[f_\varepsilon]_{res} := f_\varepsilon - [f_\varepsilon]_{ess}$$

for any measurable function $f_\varepsilon$ in $Q$.

The following estimates are nowadays standard and can be derived from the energy inequality (2.5), see [10]:

$$\left\| \left[ \frac{\varrho_\varepsilon - \overline{\varrho}}{\varepsilon} \right]_{ess} \right\|_{L^\infty(0,T,L^2(\Omega_t))} \leq C, \tag{3.1}$$

$$\|[\varrho_\varepsilon]_{res}\|_{L^\infty(0,T,L^\gamma(\Omega_t))} \leq C\varepsilon^{\frac{2}{\gamma}}, \tag{3.2}$$

$$\|[1]_{res}\|_{L^\infty(0,T,L^1(\Omega_t))} \leq C\varepsilon^2. \tag{3.3}$$

Furthermore, from (3.2) and (3.3) we derive

$$\left\| \left[ \frac{\varrho_\varepsilon - \overline{\varrho}}{\varepsilon} \right]_{res} \right\|_{L^\infty(0,T,L^q(\Omega_t))} \leq C\varepsilon^{\frac{2-q}{q}} \tag{3.4}$$

for any $1 \leq q \leq \min\{\gamma, 2\}$. Clearly, the relations (3.1), (3.4) imply (2.7).

Finally, using a version of Korn's inequality we obtain

$$\|\mathbf{u}_\varepsilon\|_{L^2(0,T,W^{1,2}(\Omega_t))} \leq C,$$

$$\|\sqrt{\varrho_\varepsilon}\mathbf{u}_\varepsilon\|_{L^\infty(0,T,L^2(\Omega_t))} \leq C,$$

where the constant $C$ does not depend on $\varepsilon$.

## 4. Weak convergence

It is convenient to prolong the quantities defined on $\Omega_t$ to the whole space $\mathbb{R}^3$. Therefore we define $\varrho := \overline{\varrho}$ on $\mathbb{R}^3 \setminus \Omega_t$ for all $t \in [0, T]$. To extend the velocity we use the standard extension operator $E_t : W^{1,2}(\Omega_t) \mapsto W^{1,2}(\mathbb{R}^3)$ which is uniformly bounded with respect to $t \in [0, T]$. With this convention we conclude from the uniform estimates in Section 3 that

$$\varrho_\varepsilon \to \overline{\varrho} \quad \text{in } L^\infty(0, T, L^\gamma(\mathbb{R}^3)) \tag{4.1}$$

$$\varrho_\varepsilon \to \overline{\varrho} \quad \text{in } C(0, T, L^r(\mathbb{R}^3)), \quad r \in [1, \gamma) \tag{4.2}$$

$$E_t\mathbf{u}_\varepsilon \to \mathbf{U} \quad \text{weakly in } L^2(0, T, W^{1,2}(\mathbb{R}^3)). \tag{4.3}$$

Moreover we get for any set $[T_1, T_2] \times K \subset Q$ with compact $K$

$$\varrho_\varepsilon\mathbf{u}_\varepsilon \to \overline{\varrho}\mathbf{U} \quad \text{weakly}^* \text{ in } L^\infty(T_1, T_2, L^{\frac{2\gamma}{\gamma+1}}(K)). \tag{4.4}$$

This allows us to proceed to the limit with all terms in (2.4) except for the convective term, for which we only have

$$\varrho_\varepsilon\mathbf{u}_\varepsilon \otimes \mathbf{u}_\varepsilon \to \overline{\varrho\mathbf{u} \otimes \mathbf{u}} \quad \text{weakly in } L^q(T_1, T_2, L^q(K)), \tag{4.5}$$

for a certain $q > 1$. Here we use the assumption $\gamma > \frac{3}{2}$ and the notation $\overline{f(v)}$ for a weak limit of the sequence $f(v_\varepsilon)$.

In the rest of this paper we prove that

$$\mathbf{u}_\varepsilon \to \mathbf{U} \quad \text{strongly in } L^2(K), \tag{4.6}$$

for every $K \subset Q$ compact. This immediately implies the desired convergence of the convective term and finishes the proof of Theorem 2.5.

## 5. Lighthill acoustic analogy

We reformulate the system (1.7)-(1.8) in the form of the Lighthill acoustic analogy. To this end we define

$$r_\varepsilon := \frac{\varrho_\varepsilon - \overline{\varrho}}{\varepsilon} \tag{5.1}$$

$$\mathbf{V}_\varepsilon := \varrho_\varepsilon \mathbf{u}_\varepsilon - \overline{\varrho}\mathbf{V} \tag{5.2}$$

and these quantities fulfill

$$\varepsilon \partial_t r_\varepsilon + \mathrm{div}_x \mathbf{V}_\varepsilon = 0 \tag{5.3}$$

$$\varepsilon \partial_t \mathbf{V}_\varepsilon + p'(\overline{\varrho})\nabla_x r_\varepsilon = \varepsilon(\mathrm{div}_x \mathbb{F}_1 + \mathbf{F}_2 + \nabla_x F_3), \tag{5.4}$$

where

$$\mathbb{F}_1 = \mathbb{S}(\nabla_x \mathbf{u}_\varepsilon) - \varrho_\varepsilon \mathbf{u}_\varepsilon \otimes \mathbf{u}_\varepsilon \tag{5.5}$$

$$\mathbf{F}_2 = -\overline{\varrho}\partial_t \mathbf{V} \tag{5.6}$$

$$F_3 = \frac{1}{\varepsilon^2}\left(p(\varrho_\varepsilon) - p'(\overline{\varrho})(\varrho_\varepsilon - \overline{\varrho}) - p(\overline{\varrho})\right). \tag{5.7}$$

Next, we make the change of coordinates $y = x - m(t)$ and rewrite the system (5.3)–(5.4) from the time dependent domain $Q = (0,T) \times \Omega_t$ to the fixed domain $\widetilde{Q} = (0,T) \times \Omega_0$. Denoting $\widetilde{f}(t,y) = f(t, y + m(t))$ for any scalar, vector or tensor-valued quantity $f$, the Lighthill acoustic analogy (cf. Lighthill [12], [13]) takes a form

$$\varepsilon \partial_t \widetilde{r}_\varepsilon + \mathrm{div}_y(\widetilde{\mathbf{V}}_\varepsilon - m'(t)\varepsilon \widetilde{r}_\varepsilon) \;=\; 0 \tag{5.8}$$

$$\varepsilon \partial_t \widetilde{\mathbf{V}}_\varepsilon + p'(\overline{\rho})\nabla_y \widetilde{r}_\varepsilon - \varepsilon\, \mathrm{div}_y(\widetilde{\mathbf{V}}_\varepsilon \otimes m'(t)) \;=\; \varepsilon\left(\mathrm{div}_y \widetilde{\mathbb{F}}_1 + \widetilde{\mathbf{F}}_2 + \nabla_y \widetilde{F}_3\right). \tag{5.9}$$

Morevorer, we denote $\widetilde{\mathbf{W}}_\varepsilon = (\widetilde{\mathbf{V}}_\varepsilon - m'(t)\varepsilon\widetilde{r}_\varepsilon)$. Then we rewrite (5.8)–(5.9) further to

$$\varepsilon \partial_t \widetilde{r}_\varepsilon + \mathrm{div}_y \widetilde{\mathbf{W}}_\varepsilon = 0 \tag{5.10}$$

$$\varepsilon \partial_t \widetilde{\mathbf{W}}_\varepsilon + p'(\overline{\rho})\nabla_y \widetilde{r}_\varepsilon = \varepsilon\left(\mathrm{div}_y(\widetilde{\mathbb{F}}_1 + \widetilde{\mathbf{V}}_\varepsilon \otimes m'(t)) + \widetilde{\mathbf{F}}_2 + \nabla_y \widetilde{F}_3 - \varepsilon m''(t)\widetilde{r}_\varepsilon + \varepsilon m'(t)\, \mathrm{div}_y \widetilde{\mathbf{W}}_\varepsilon\right). \tag{5.11}$$

Note in particular that such defined $\widetilde{\mathbf{W}}_\varepsilon$ satisfies also the boundary condition $\widetilde{\mathbf{W}}_\varepsilon \cdot \mathbf{n} = 0$ on $\partial\Omega_0$. The weak formulation of (5.10)–(5.11) reads

$$\int_0^T \int_{\Omega_0} \varepsilon \widetilde{r}_\varepsilon \partial_t \varphi + \widetilde{\mathbf{W}}_\varepsilon \cdot \nabla_y \varphi \; \mathrm{d}y\, \mathrm{d}t = -\int_{\Omega_0} \varepsilon \varrho^{(1)}_{0,\varepsilon} \varphi(0,\cdot)\; \mathrm{d}y \tag{5.12}$$

for all $\varphi \in C_c^\infty([0,T) \times \overline{\Omega}_0)$ and

$$\begin{aligned}
&\int_0^T \int_{\Omega_0} \varepsilon \widetilde{\mathbf{W}}_\varepsilon \cdot \partial_t \varphi + p'(\overline{\varrho}) r_\varepsilon \operatorname{div}_y \varphi \, \mathrm{d}y \, \mathrm{d}t \\
&= - \int_{\Omega_0} \Big( \varepsilon \overline{\varrho} (\mathbf{u}_{0,\varepsilon} - \mathbf{V}(0,\cdot)) + \varepsilon^2 \varrho_{0,\varepsilon}^{(1)} (\mathbf{u}_{0,\varepsilon} - m'(0)) \Big) \varphi(0,\cdot) \, \mathrm{d}y \\
&+ \varepsilon \int_0^T \int_{\Omega_0} (\widetilde{\mathbb{F}}_1 + \widetilde{\mathbf{V}}_\varepsilon \otimes m'(t) + \varepsilon m'(t) \otimes \widetilde{\mathbf{W}}_\varepsilon) : \nabla_y \varphi - (\widetilde{\mathbf{F}}_2 - \varepsilon m''(t) \widetilde{r}_\varepsilon) \cdot \varphi + \widetilde{F}_3 \operatorname{div}_y \varphi \, \mathrm{d}y \, \mathrm{d}t
\end{aligned} \tag{5.13}$$

for all $\varphi \in C_c^\infty([0,T) \times \overline{\Omega}_0)$ such that $\varphi \cdot \mathbf{n} = 0$ on $\partial \Omega_0$.

5.1. **Helmholtz projection and Neumann Laplacian.** For any $\mathbf{v} \in L^p(\Omega_0)$ we denote by $\mathbf{H}(\mathbf{v})$ its Helmholtz projection, more precisely

$$\mathbf{H}(\mathbf{v}) := \mathbf{v} - \nabla_y \Theta, \tag{5.14}$$

where $\Theta$ such that $\nabla_y \Theta \in L^p(\Omega_0)$ is a unique solution of the problem

$$\Delta \Theta = \operatorname{div}_y \mathbf{v}, \quad \frac{\partial \Theta}{\partial \mathbf{n}} = \mathbf{v} \cdot \mathbf{n} \text{ on } \partial \Omega_0, \quad |\Theta| \to 0 \text{ as } |y| \to \infty \tag{5.15}$$

which in weak formulation reads

$$\int_{\Omega_0} \nabla_y \Theta \cdot \nabla_y \varphi \, \mathrm{d}y = \int_{\Omega_0} \mathbf{v} \cdot \nabla_y \varphi \, \mathrm{d}y \quad \text{for all } \varphi \in C_c^\infty(\overline{\Omega}_0). \tag{5.16}$$

Neumann Laplacian operator $\Delta_N$ plays a crucial role in the following analysis. We recall that $-\Delta_N$ is a nonnegative self-adjoint operator on $L^2(\Omega_0)$ with domain

$$\mathcal{D}(-\Delta_N) = \left\{ w \in W^{2,2}(\Omega_0), \nabla_y w \cdot \mathbf{n} = 0 \text{ on } \partial \Omega_0 \right\}. \tag{5.17}$$

5.2. **Compactness of the solenoidal part.** Our aim is to prove that for $\varphi \in C_c^\infty([0,T] \times \overline{\Omega}_0)$ with $\varphi \cdot \mathbf{n} = 0$ on $\partial \Omega_0$ it holds,

$$\left[ t \mapsto \int_{\Omega_0} \widetilde{\mathbf{W}}_\varepsilon \cdot \varphi \, \mathrm{d}y \right] \longrightarrow \left[ t \mapsto \int_{\Omega_0} \widetilde{\mathbf{W}} \cdot \varphi \, \mathrm{d}y \right] \text{ strongly in } L^2(0,T), \tag{5.18}$$

with $\widetilde{\mathbf{W}} := \overline{\varrho}(\widetilde{\mathbf{U}} - \widetilde{\mathbf{V}})$.

Using the Helmholtz decomposition in $\Omega_0$ we split $\varphi = \mathbf{H}(\varphi) + \mathbf{H}^\perp(\varphi)$. Due to estimates from section 3 we have that $\int_{\Omega_0} \widetilde{\mathbf{W}}_\varepsilon(t,\cdot) \cdot \mathbf{H}(\varphi) \, \mathrm{d}y$ is bounded independently of $t$ and $\varepsilon$. Further, due to (5.11) and an Aubin-Lions argument we conclude

$$\left[ t \mapsto \int_{\Omega_0} \widetilde{\mathbf{W}}_\varepsilon \cdot \mathbf{H}(\varphi) \, \mathrm{d}y \right] \longrightarrow \left[ t \mapsto \int_{\Omega_0} \widetilde{\mathbf{W}} \cdot \varphi \, \mathrm{d}y \right] \text{ strongly in } L^2(0,T), \tag{5.19}$$

In the rest of this paper we discuss the gradient part of the velocity.

## 6. Convergence of the gradient part of velocity

We introduce the acoustic potential $\Psi_\varepsilon$ as the gradient part of the quantity $\widetilde{\mathbf{W}}_\varepsilon$, more precisely

$$\widetilde{\mathbf{W}}_\varepsilon = \mathbf{H}(\widetilde{\mathbf{W}}_\varepsilon) + \nabla_y \Psi_\varepsilon. \tag{6.1}$$

Observe that $\nabla_y \Delta_N^{-1} \varphi$ with $\varphi \in C_c^\infty([0,T) \times \Omega_0)$ is an admissible test function in equation (5.13). Using this test function and having in mind the relation

$$\int_{\Omega_0} \nabla_y \Psi_\varepsilon \cdot \nabla_y \varphi \, \mathrm{d}y = \int_{\Omega_0} \widetilde{\mathbf{W}}_\varepsilon \cdot \nabla_y \varphi \, \mathrm{d}y \quad \text{ for all } \varphi \in C_c^\infty(\overline{\Omega}_0) \tag{6.2}$$

we obtain from (5.12)–(5.13)

$$\int_0^T \int_{\Omega_0} \varepsilon \widetilde{r}_\varepsilon \partial_t \varphi + \nabla_y \Psi_\varepsilon \cdot \nabla_y \varphi \, \mathrm{d}y \, \mathrm{d}t = -\int_{\Omega_0} \varepsilon \varrho_{0,\varepsilon}^{(1)} \varphi(0,\cdot) \, \mathrm{d}y \tag{6.3}$$

for all $\varphi \in C_c^\infty([0,T) \times \overline{\Omega}_0)$ and

$$\begin{aligned}
&\int_0^T \int_{\Omega_0} \varepsilon \Psi_\varepsilon \partial_t \varphi + p'(\overline{\varrho}) r_\varepsilon \varphi \, \mathrm{d}y \, \mathrm{d}t \\
= &-\int_{\Omega_0} \left( \varepsilon \overline{\varrho} (\mathbf{u}_{0,\varepsilon} - \mathbf{V}(0,\cdot)) + \varepsilon^2 \varrho_{0,\varepsilon}^{(1)} (\mathbf{u}_{0,\varepsilon} - m'(0)) \right) \cdot \nabla_y \Delta_N^{-1} \varphi(0,\cdot) \, \mathrm{d}y \\
&+\varepsilon \int_0^T \int_{\Omega_0} (\widetilde{\mathbb{F}}_1 + \widetilde{\mathbf{V}}_\varepsilon \otimes m'(t) + \varepsilon m'(t) \otimes \widetilde{\mathbf{W}}_\varepsilon) : \nabla_y^2 \Delta_N^{-1} \varphi \, \mathrm{d}y \, \mathrm{d}t \\
&-\varepsilon \int_0^T \int_{\Omega_0} (\widetilde{\mathbf{F}}_2 - \varepsilon m''(t) \widetilde{r}_\varepsilon) \cdot \nabla_y \Delta_N^{-1} \varphi + \widetilde{F}_3 \varphi \, \mathrm{d}y \, \mathrm{d}t =: \varepsilon [h_\varepsilon, \varphi]
\end{aligned} \tag{6.4}$$

for all $\varphi \in C_c^\infty([0,T) \times \overline{\Omega}_0)$.

6.1. **Uniform bounds revisited.** It is easy to observe that all uniform bounds from Section 3 transfer from the time dependent domain $(0,T) \times \Omega_t$ to fixed domain $(0,T) \times \Omega_0$. Moreover we deduce easily

$$\| [\widetilde{r}_\varepsilon]_{ess} \|_{L^\infty(0,T,L^2(\Omega_0))} \leq C \tag{6.5}$$

$$\| [\widetilde{r}_\varepsilon]_{res} \|_{L^\infty(0,T,L^q(\Omega_0))} \leq C \varepsilon^{\frac{2-q}{q}} \tag{6.6}$$

for any $1 \leq q < \min\{\gamma, 2\}$. Moreover we have

$$\left\| [\sqrt{\widetilde{\varrho}_\varepsilon}]_{ess} \sqrt{\widetilde{\varrho}_\varepsilon} \widetilde{\mathbf{u}}_\varepsilon \right\|_{L^\infty(0,T,L^2(\Omega_0))} \leq C \tag{6.7}$$

$$\left\| [\sqrt{\widetilde{\varrho}_\varepsilon}]_{res} \sqrt{\widetilde{\varrho}_\varepsilon} \widetilde{\mathbf{u}}_\varepsilon \right\|_{L^\infty(0,T,L^q(\Omega_0))} \leq C \varepsilon^{\frac{1}{\gamma}} \tag{6.8}$$

for $q = \frac{2\gamma}{\gamma+1}$.

6.2. **Estimate of forcing term.** First, we estimate the terms $\widetilde{\mathbb{F}}_1$ and $\widetilde{F}_3$ in the same manner as in [7, Section 4]. Observing that

$$\begin{aligned}
&\left|\int_{\Omega_0} \mathbb{S}(\nabla_y \widetilde{\mathbf{u}}_\varepsilon) : \nabla_y^2 \Delta_N^{-1} \varphi \, \mathrm{d}y\right| \\
\leq & C \left\| \mathbb{S}(\nabla_y \widetilde{\mathbf{u}}_\varepsilon) \right\|_{L^2(\Omega_0)} \left\| \nabla_y^2 \Delta_N^{-1} \varphi \right\|_{L^2(\Omega_0)} \\
\leq & C \left\| \mathbb{S}(\nabla_y \widetilde{\mathbf{u}}_\varepsilon) \right\|_{L^2(\Omega_0)} \left( \|\varphi\|_{L^2(\Omega_0)} + \left\| (-\Delta_N)^{-1} \varphi \right\|_{L^2(\Omega_0)} \right),
\end{aligned} \tag{6.9}$$

the Riesz representation theorem yields the existence of functions $F_{i,\varepsilon} \in L^2((0,T) \times \Omega_0)$ for $i = 1,2$ such that

$$\int_0^T \int_{\Omega_0} \mathbb{S}(\nabla_y \widetilde{\mathbf{u}}_\varepsilon) : \nabla_y^2 \Delta_N^{-1} \varphi \, \mathrm{d}y \, \mathrm{d}t = \int_0^T \int_{\Omega_0} F_{1,\varepsilon} \varphi + F_{2,\varepsilon} (-\Delta_N)^{-1} \varphi \, \mathrm{d}y \, \mathrm{d}t. \tag{6.10}$$

Similarly we proceed with the convective term. Here we write

$$\begin{aligned}
&\int_{\Omega_0} \widetilde{\varrho}_\varepsilon \widetilde{\mathbf{u}}_\varepsilon \otimes \widetilde{\mathbf{u}}_\varepsilon : \nabla_y^2 \Delta_N^{-1} \varphi \, \mathrm{d}y \\
= & \int_{\Omega_0} [\widetilde{\varrho}_\varepsilon]_{ess} \widetilde{\mathbf{u}}_\varepsilon \otimes \widetilde{\mathbf{u}}_\varepsilon : \nabla_y^2 \Delta_N^{-1} \varphi \, \mathrm{d}y + \int_{\Omega_0} [\sqrt{\widetilde{\varrho}_\varepsilon}]_{res} \sqrt{\widetilde{\varrho}_\varepsilon} \widetilde{\mathbf{u}}_\varepsilon \otimes \widetilde{\mathbf{u}}_\varepsilon : \nabla_y^2 \Delta_N^{-1} \varphi \, \mathrm{d}y.
\end{aligned} \tag{6.11}$$

We estimate the essential part as follows

$$\begin{aligned}
&\left|\int_{\Omega_0} [\widetilde{\varrho}_\varepsilon]_{ess} \widetilde{\mathbf{u}}_\varepsilon \otimes \widetilde{\mathbf{u}}_\varepsilon : \nabla_y^2 \Delta_N^{-1} \varphi \, \mathrm{d}y\right| \\
\leq & \left\| [\widetilde{\varrho}_\varepsilon]_{ess} \widetilde{\mathbf{u}}_\varepsilon \right\|_{L^2(\Omega_0)} \left\| \widetilde{\mathbf{u}}_\varepsilon \right\|_{L^6(\Omega_0)} \left\| \nabla_y^2 \Delta_N^{-1} \varphi \right\|_{L^3(\Omega_0)} \\
\leq & \left\| [\widetilde{\varrho}_\varepsilon]_{ess} \widetilde{\mathbf{u}}_\varepsilon \right\|_{L^2(\Omega_0)} \left\| \widetilde{\mathbf{u}}_\varepsilon \right\|_{L^6(\Omega_0)} \left( \|\varphi\|_{L^3(\Omega_0)} + \left\| (-\Delta_N)^{-1} \varphi \right\|_{L^3(\Omega_0)} \right),
\end{aligned} \tag{6.12}$$

Using the interpolation inequality we easily have

$$\|\varphi\|_{L^3(\Omega_0)} \leq C \left( \|\varphi\|_{L^2(\Omega_0)} + \|\varphi\|_{L^6(\Omega_0)} \right) \leq C \left( \|\varphi\|_{L^2(\Omega_0)} + \left\| (-\Delta_N)^{1/2} \varphi \right\|_{L^2(\Omega_0)} \right) \tag{6.13}$$

and similarly

$$\left\| (-\Delta_N)^{-1} \varphi \right\|_{L^3(\Omega_0)} \leq C \left( \left\| (-\Delta_N)^{-1} \varphi \right\|_{L^2(\Omega_0)} + \left\| (-\Delta_N)^{-1/2} \varphi \right\|_{L^2(\Omega_0)} \right) \tag{6.14}$$

and thus there exist functions $F_{i,\varepsilon} \in L^2((0,T) \times \Omega_0)$ for $i = 3, ..., 6$ such that

$$\begin{aligned}
&\int_0^T \int_{\Omega_0} [\widetilde{\varrho}_\varepsilon]_{ess} \widetilde{\mathbf{u}}_\varepsilon \otimes \widetilde{\mathbf{u}}_\varepsilon : \nabla_y^2 \Delta_N^{-1} \varphi \, \mathrm{d}y \\
= & \int_0^T \int_{\Omega_0} \sum_{i=3}^{6} F_{i,\varepsilon} (-\Delta_N)^{2-i/2} \varphi \, \mathrm{d}y \, \mathrm{d}t.
\end{aligned} \tag{6.15}$$

For the residual part of the convective term we proceed as follows

$$\begin{aligned}
&\left|\int_{\Omega_0} [\sqrt{\widetilde{\varrho}_\varepsilon}]_{res}\sqrt{\widetilde{\varrho}_\varepsilon}\widetilde{\mathbf{u}}_\varepsilon \otimes \widetilde{\mathbf{u}}_\varepsilon : \nabla_y^2 \Delta_N^{-1}\varphi \ \mathrm{d}y\right| \\
&\leq \left\|[\sqrt{\widetilde{\varrho}_\varepsilon}]_{res}\sqrt{\widetilde{\varrho}_\varepsilon}\widetilde{\mathbf{u}}_\varepsilon\right\|_{L^q(\Omega_0)} \|\widetilde{\mathbf{u}}_\varepsilon\|_{L^6(\Omega_0)} \left\|\nabla_y^2\Delta_N^{-1}\varphi\right\|_{L^r(\Omega_0)} \\
&\leq \left\|[\sqrt{\widetilde{\varrho}_\varepsilon}]_{res}\sqrt{\widetilde{\varrho}_\varepsilon}\widetilde{\mathbf{u}}_\varepsilon\right\|_{L^q(\Omega_0)} \|\widetilde{\mathbf{u}}_\varepsilon\|_{L^6(\Omega_0)} \left(\|\varphi\|_{L^r(\Omega_0)} + \left\|(-\Delta_N)^{-1}\varphi\right\|_{L^r(\Omega_0)}\right), \qquad (6.16)
\end{aligned}$$

where $q = \frac{2\gamma}{\gamma+1}$ and

$$\frac{\gamma+1}{2\gamma} + \frac{1}{6} + \frac{1}{r} = 1. \qquad (6.17)$$

Note that $r > 3$, so we estimate the arising norms of the test function

$$\|\varphi\|_{L^r(\Omega_0)} \leq C\left(\|\varphi\|_{L^2(\Omega_0)} + \left\|\nabla_y^2\varphi\right\|_{L^2(\Omega_0)}\right) \leq C\left(\|\varphi\|_{L^2(\Omega_0)} + \|(-\Delta_N)\varphi\|_{L^2(\Omega_0)}\right) \qquad (6.18)$$

and similarly

$$\left\|(-\Delta_N)^{-1}\varphi\right\|_{L^r(\Omega_0)} \leq C\left(\left\|(-\Delta_N)^{-1}\varphi\right\|_{L^2(\Omega_0)} + \|\varphi\|_{L^2(\Omega_0)}\right). \qquad (6.19)$$

This again yields the existence of functions $F_{i,\varepsilon} \in L^2((0,T)\times\Omega_0)$ for $i = 7, ..., 9$ such that

$$\begin{aligned}
&\int_0^T \int_{\Omega_0} [\sqrt{\widetilde{\varrho}_\varepsilon}]_{res}\sqrt{\widetilde{\varrho}_\varepsilon}\widetilde{\mathbf{u}}_\varepsilon \otimes \widetilde{\mathbf{u}}_\varepsilon : \nabla_y^2\Delta_N^{-1}\varphi \ \mathrm{d}y\,\mathrm{d}t \\
&= \int_0^T \int_{\Omega_0} F_{7,\varepsilon}(-\Delta_N)\varphi + F_{8,\varepsilon}\varphi + F_{9,\varepsilon}(-\Delta_N)^{-1}\varphi \ \mathrm{d}y\,\mathrm{d}t. \qquad (6.20)
\end{aligned}$$

Next, we estimate in a similar manner the pressure term $\widetilde{F}_3$. Since it holds

$$\begin{aligned}
&\left|\int_{\Omega_0} \frac{1}{\varepsilon^2}\left(p(\widetilde{\varrho}_\varepsilon) - p'(\overline{\varrho})(\widetilde{\varrho}_\varepsilon - \overline{\varrho}) - p(\overline{\varrho})\right)\varphi \ \mathrm{d}y\right| \\
&\leq \left\|\frac{1}{\varepsilon^2}\left(p(\widetilde{\varrho}_\varepsilon) - p'(\overline{\varrho})(\widetilde{\varrho}_\varepsilon - \overline{\varrho}) - p(\overline{\varrho})\right)\right\|_{L^1(\Omega_0)} \|\varphi\|_{L^\infty(\Omega_0)} \qquad (6.21)
\end{aligned}$$

we have to estimate the $L^\infty$ norm of $\varphi$ in terms of $L^2$ norms of powers of $(-\Delta_N)\varphi$. We have

$$\|\varphi\|_{L^\infty(\Omega_0)} \leq C(\|\nabla_y\varphi\|_{L^6(\Omega_0)} + \|\varphi\|_{L^6(\Omega_0)}) \leq C(\left\|\nabla_y^2\varphi\right\|_{L^2(\Omega_0)} + \|\nabla_y\varphi\|_{L^2(\Omega_0)}) \qquad (6.22)$$

and again we use the inequality

$$\left\|\nabla_y^2\varphi\right\|_{L^2(\Omega_0)} \leq C(\|\varphi\|_{L^2(\Omega_0)} + \|(-\Delta_N)\varphi\|_{L^2(\Omega_0)}) \qquad (6.23)$$

together with

$$\|\nabla_y\varphi\|_{L^2(\Omega_0)} = \left\|(-\Delta_N)^{1/2}\varphi\right\|_{L^2(\Omega_0)} \qquad (6.24)$$

to conclude that there exist functions $F_{i,\varepsilon} \in L^2((0,T)\times\Omega_0)$ for $i = 10, ..., 12$ such that

$$\begin{aligned}
&\int_0^T \int_{\Omega_0} \frac{1}{\varepsilon^2}\left(p(\widetilde{\varrho}_\varepsilon) - p'(\overline{\varrho})(\widetilde{\varrho}_\varepsilon - \overline{\varrho}) - p(\overline{\varrho})\right)\varphi \ \mathrm{d}y\,\mathrm{d}t \\
&= \int_0^T \int_{\Omega_0} F_{10,\varepsilon}(-\Delta_N)\varphi + F_{11,\varepsilon}(-\Delta_N)^{1/2}\varphi + F_{12,\varepsilon}\varphi \ \mathrm{d}y\,\mathrm{d}t. \qquad (6.25)
\end{aligned}$$

It is easy to estimate the term $\widetilde{\mathbf{F}}_2$ as follows

$$\left|\int_{\Omega_0} \overline{\varrho}\partial_t \widetilde{\mathbf{V}} \cdot \nabla_y \Delta_N^{-1}\varphi \, \mathrm{d}y\right| \leq C \left\|\partial_t \widetilde{\mathbf{V}}\right\|_{L^2(\Omega_0)} \left\|\nabla_y \Delta_N^{-1}\varphi\right\|_{L^2(\Omega_0)} \tag{6.26}$$

and

$$\left\|\nabla_y \Delta_N^{-1}\varphi\right\|_{L^2(\Omega_0)} = \left\|(-\Delta_N)^{1/2}\varphi\right\|_{L^2(\Omega_0)}, \tag{6.27}$$

thus there exists a function $F_{13} \in L^2((0,T \times \Omega_0))$ such that

$$\int_0^T \int_{\Omega_0} \overline{\varrho}\partial_t \widetilde{\mathbf{V}} \cdot \nabla_y \Delta_N^{-1}\varphi \, \mathrm{d}y \, \mathrm{d}t = \int_0^T \int_{\Omega_0} F_{13}(-\Delta_N)^{1/2}\varphi \, \mathrm{d}y \, \mathrm{d}t. \tag{6.28}$$

Summing up, we proved up to now the existence of functions $G_{\varepsilon,i} \in L^2(0,T,L^2(\Omega_0))$, $i \in \{1,\dots,5\}$ such that

$$\int_0^T \int_{\Omega_0} \widetilde{\mathbb{F}}_1 : \nabla_y^2 \Delta_N^{-1}\varphi - \widetilde{\mathbf{F}}_2 \cdot \nabla_y \Delta_N^{-1}\varphi + \widetilde{F}_3\varphi \, \mathrm{d}y \, \mathrm{d}t = \sum_{i=1}^5 \int_0^T \int_{\Omega_0} G_{\varepsilon,i}(-\Delta_N)^{-1+\frac{i-1}{2}}[\varphi]. \tag{6.29}$$

Moreover, there exists $c \in \mathbb{R}$ independent of $\varepsilon$ and $i$ fulfilling

$$\sum_{i=1}^5 \|G_{\varepsilon,i}\|_{L^2(0,T,L^2(\Omega_0))} \leq c. \tag{6.30}$$

Now we estimate the extra terms due to the translation of the domain. Again we first split $\widetilde{\mathbf{V}}_\varepsilon$ to the essential and residual part

$$\widetilde{\mathbf{V}}_\varepsilon = [\widetilde{\varrho}_\varepsilon]_{ess}\widetilde{\mathbf{u}}_\varepsilon + [\widetilde{\varrho}_\varepsilon]_{res}\widetilde{\mathbf{u}}_\varepsilon - \overline{\varrho}\widetilde{\mathbf{V}},$$

and remind that

$$\|[\sqrt{\widetilde{\varrho}_\varepsilon}]_{res}\sqrt{\widetilde{\varrho}_\varepsilon}\widetilde{\mathbf{u}}_\varepsilon\|_{L^\infty(0,T,L^q(\Omega_t))} \leq c$$

for $q = \frac{2\gamma}{\gamma+1}$ and

$$\|[\sqrt{\widetilde{\varrho}_\varepsilon}]_{ess}\sqrt{\widetilde{\varrho}_\varepsilon}\widetilde{\mathbf{u}}_\varepsilon\|_{L^\infty(0,T,L^2(\Omega_t))} \leq c.$$

Consequently, we estimate the term $\widetilde{\mathbf{V}}_\varepsilon \otimes m'(t)$ as follows. The essential part is treated easily

$$\left|\int_{\Omega_0} \sqrt{\widetilde{\varrho}_\varepsilon}]_{ess}\sqrt{\widetilde{\varrho}_\varepsilon}\widetilde{\mathbf{u}}_\varepsilon \otimes m'(t) : \nabla_y^2 \Delta_N^{-1}\varphi \, \mathrm{d}y\right| \leq C \left\|\sqrt{\widetilde{\varrho}_\varepsilon}]_{ess}\sqrt{\widetilde{\varrho}_\varepsilon}\widetilde{\mathbf{u}}_\varepsilon\right\|_{L^2(\Omega_0)} \left\|\nabla_y^2 \Delta_N^{-1}\varphi\right\|_{L^2(\Omega_0)}$$
$$\leq \left\|\sqrt{\widetilde{\varrho}_\varepsilon}]_{ess}\sqrt{\widetilde{\varrho}_\varepsilon}\widetilde{\mathbf{u}}_\varepsilon\right\|_{L^2(\Omega_0)} \left(\|\varphi\|_{L^2(\Omega_0)} + \left\|(-\Delta_N)^{-1}\varphi\right\|_{L^2(\Omega_0)}\right) \tag{6.31}$$

and the term $\overline{\varrho}\widetilde{\mathbf{V}} \otimes m'(t)$ is estimated in the same way, because it also belongs to the space $L^2(0,T,L^2(\Omega_0))$. The residual part is estimated similarly as in the case of the convective term.

$$\left|\int_{\Omega_0} \sqrt{\widetilde{\varrho}_\varepsilon}]_{res}\sqrt{\widetilde{\varrho}_\varepsilon}\widetilde{\mathbf{u}}_\varepsilon \otimes m'(t) : \nabla_y^2 \Delta_N^{-1}\varphi \, \mathrm{d}y\right| \leq C \left\|\sqrt{\widetilde{\varrho}_\varepsilon}]_{res}\sqrt{\widetilde{\varrho}_\varepsilon}\widetilde{\mathbf{u}}_\varepsilon\right\|_{L^q(\Omega_0)} \left\|\nabla_y^2 \Delta_N^{-1}\varphi\right\|_{L^{q'}(\Omega_0)}$$
$$\leq \left\|\sqrt{\widetilde{\varrho}_\varepsilon}]_{res}\sqrt{\widetilde{\varrho}_\varepsilon}\widetilde{\mathbf{u}}_\varepsilon\right\|_{L^q(\Omega_0)} \left(\|\varphi\|_{L^{q'}(\Omega_0)} + \left\|(-\Delta_N)^{-1}\varphi\right\|_{L^{q'}(\Omega_0)}\right) \tag{6.32}$$

with $1/q + 1/q' = 1$. Since $\gamma > 3/2$, $q' < 6$ therefore we estimate the arising norms of the test function in the same manner as in (6.18) and (6.19).

Next term we want to estimate is $\varepsilon m'(t) \otimes \widetilde{\mathbf{W}}_\varepsilon$. Since $\widetilde{\mathbf{W}}_\varepsilon = \widetilde{\mathbf{V}}_\varepsilon - m'(t)\varepsilon\widetilde{r}_\varepsilon$, we use the estimates (6.31)–(6.32) to handle the part including $\widetilde{\mathbf{V}}_\varepsilon$. Therefore we have to estimate $\varepsilon\widetilde{r}_\varepsilon$. Again we split $\widetilde{r}_\varepsilon$ to the essential part and the residual part

$$\widetilde{r}_\varepsilon = [\widetilde{r}_\varepsilon]_{ess} + [\widetilde{r}_\varepsilon]_{res},$$

and estimate them separately. The essential part belongs to $L^\infty(0,T,L^2(\Omega_0))$ and its estimate is straightforward, whereas for the residual part we use the estimate (6.6) and estimate it in a similar manner as in (6.32).

It remains to estimate the term $m''(t)\varepsilon\widetilde{r}_\varepsilon$. Again splitting it to essential and residual part we have

$$\left| \int_{\Omega_0} m''(t)[\widetilde{r}_\varepsilon]_{ess} \cdot \nabla_y \Delta_N^{-1} \varphi \,\mathrm{d}y \right| \leq C \left\| [\widetilde{r}_\varepsilon]_{ess} \right\|_{L^2(\Omega_0)} \left\| (-\Delta_N)^{1/2} \varphi \right\|_{L^2(\Omega_0)}, \tag{6.33}$$

where we used (6.24), whereas

$$\begin{aligned} &\left| \int_{\Omega_0} m''(t)[\widetilde{r}_\varepsilon]_{res} \cdot \nabla_y \Delta_N^{-1} \varphi \,\mathrm{d}y \right| \leq C \left\| [\widetilde{r}_\varepsilon]_{ess} \right\|_{L^{\frac{6}{5}}(\Omega_0)} \left\| \nabla_y (-\Delta_N) \varphi \right\|_{L^6(\Omega_0)} \\ \leq &C \left\| [\widetilde{r}_\varepsilon]_{ess} \right\|_{L^{\frac{6}{5}}(\Omega_0)} \left\| \nabla_y^2 (-\Delta_N) \varphi \right\|_{L^2(\Omega_0)} \leq C \left\| [\widetilde{r}_\varepsilon]_{ess} \right\|_{L^{\frac{6}{5}}(\Omega_0)} \left( \| (-\Delta_N) \varphi \|_{L^2(\Omega_0)} + \| \varphi \|_{L^2(\Omega_0)} \right). \end{aligned} \tag{6.34}$$

Therefore we can finally write

$$[h_\varepsilon, \varphi] = \sum_{i=1}^{5} \int_0^T \int_{\Omega_0} G_{\varepsilon,i} (-\Delta_N)^{-1+\frac{i-1}{2}} \varphi \,\mathrm{d}y \,\mathrm{d}t, \tag{6.35}$$

with $G_{\varepsilon,i}$ satisfying

$$\sum_{i=1}^{5} \| G_{\varepsilon,i} \|_{L^2(0,T,L^2(\Omega_0))} \leq c. \tag{6.36}$$

.

6.3. **Solution to wave equation.** The system (6.3)–(6.4) can be understood as

$$\begin{aligned} \varepsilon \partial_t \widetilde{r}_\varepsilon - (-\Delta_N) \Psi_\varepsilon &= 0 & (6.37) \\ \varepsilon \partial_t \Psi_\varepsilon + p'(\overline{\varrho}) \widetilde{r}_\varepsilon &= \varepsilon h_\varepsilon, & (6.38) \end{aligned}$$

keeping in mind also the boundary condition

$$\nabla_y \Psi_\varepsilon \cdot \mathbf{n} = 0 \quad \text{ on } \partial\Omega_0. \tag{6.39}$$

We establish the following notation

$$w_\varepsilon = \begin{pmatrix} \widetilde{r}_\varepsilon \\ \Psi_\varepsilon \end{pmatrix}, \quad A = \begin{pmatrix} 0 & -(-\Delta_N) \\ p'(\overline{\varrho}) & 0 \end{pmatrix}, \quad b = \begin{pmatrix} 0 \\ \varepsilon h_\varepsilon \end{pmatrix}$$

and then the solution to (6.37)–(6.38) is obtained by means of the Duhamel's formula

$$w_\varepsilon(t) = e^{-\frac{1}{\varepsilon}tA} \left( \int_0^t e^{\frac{1}{\varepsilon}sA} \frac{1}{\varepsilon} b(s) \,\mathrm{d}s \right) + e^{-\frac{1}{\varepsilon}tA} w_0. \tag{6.40}$$

Hereinafter we write $e_+(t) := e^{\frac{\sqrt{p'(\overline{\varrho})}}{\varepsilon} t \mathrm{i}\sqrt{-\Delta_N}}$ and $e_-(t) := \frac{1}{e_+}$. From (6.40) we get in particular

$$\begin{aligned}\Psi_\varepsilon(t) = \frac{\sqrt{p(\overline{\varrho})}}{2\mathrm{i}\sqrt{-\Delta_N}} \left(e_-(t) - e_+(t)\right) \widetilde{r}_\varepsilon(0) + \frac{1}{2}\left(e_+(t) + e_-(t)\right) \Psi_\varepsilon(0) \\ + \frac{e_+(t) - e_-(t)}{4} \int_0^t \left(e_-(s) - e_+(s)\right) h_\varepsilon(s)\,\mathrm{d}s \\ + \frac{1}{4}\left(e_+(t) + e_-(t)\right) \int_0^t \left(e_+(s) + e_-(s)\right) h_\varepsilon(s)\,\mathrm{d}s. \end{aligned} \tag{6.41}$$

We achieve the desired strong convergence of the velocities using the RAGE theorem (see [2, Section 5.4], [17, Theorem XI.115])

**Theorem 6.1.** *Let $A$ be a self-adjoint operator and let $C$ be a bounded operator such that $C(A+\mathrm{i})^{-1}$ is compact. Then*

$$\int_0^T e^{\mathrm{i}A\frac{t}{\varepsilon}} C P_c e^{-\mathrm{i}A\frac{t}{\varepsilon}}\,\mathrm{d}t \to 0 \ \textit{as}\ \varepsilon \to 0, \tag{6.42}$$

*where $P_C$ is the projection onto the orthogonal complement of the eigenvectors of $A$.*

Further we proceed as in [6]. We apply (6.42) to $A = \sqrt{p(\overline{\varrho})}\sqrt{-\Delta_N}$, $C = \chi^2 G(-\Delta_N)$ with $\chi \in C_c^\infty(\mathbb{R}^3 \times \mathbb{R}_0^+)$, $G \in C_c^\infty(0,\infty)$, $0 \le G \le 1$. We get for $X, Y \in L^2(\Omega_0)$

$$\int_0^T \left\langle e_-(t)\chi^2 G(-\Delta_N) e_+(t) X, Y\right\rangle \,\mathrm{d}t \le \omega(\varepsilon) \|X\|_{L^2(\Omega_0)} \|Y\|_{L^2(\Omega_0)},$$

with $\omega(\varepsilon) \to 0$ as $\varepsilon \to 0$. Here $\langle \cdot, \cdot \rangle$ denotes the standard scalar product in $L^2(\Omega_0)$. For $Y = G(-\Delta_N)[X]$ it holds

$$\int_0^T \|\chi G(-\Delta_N) e_\pm(t)[X]\|_2^2 \,\mathrm{d}t \le \omega(\varepsilon)\|X\|_2^2. \tag{6.43}$$

We may apply (6.43) to the right hand side of (6.41). Thus

$$\int_0^T \left\| \chi \frac{\sqrt{p(\overline{\varrho})}}{2\mathrm{i}} \frac{G(-\Delta_N)}{\sqrt{-\Delta_N}} (e_-(t) - e_+(t)) \widetilde{r}_\varepsilon(0) \right\|_2^2 \mathrm{d}t \le \omega(\varepsilon) \|\widetilde{r}_\varepsilon(0)\|_2^2 = \omega(\varepsilon) \|\varrho_{0,\varepsilon}^{(1)}\|_2^2,$$

$$\int_0^T \left\| \frac{1}{2} \chi \frac{G(-\Delta_N)}{(-\Delta_N)^{\frac{1}{2}}} (e_+(t) + e_-(t)) (-\Delta_N)^{\frac{1}{2}} \Psi_\varepsilon(0) \right\|_2^2 \mathrm{d}t \le \omega(\varepsilon) \|\nabla \Psi_\varepsilon(0)\|_2^2$$

and, finally,

$$\begin{aligned}\int_0^T \int_0^t \left\| \frac{1}{4} (e_+(t) \pm e_-(t))(e_-(s) \pm e_+(s)) \sum_{i=1}^5 \chi(-\Delta_N)^{-1+\frac{i-1}{2}} G(-\Delta_N) G_{\varepsilon,i} \right\|_2^2 \mathrm{d}s\,\mathrm{d}t \\ \le \omega(\varepsilon) \int_0^T \sum_{i=1}^5 \|G_{\varepsilon,i}\|_2^2 \ \mathrm{d}t. \end{aligned} \tag{6.44}$$

Putting these estimates together and choosing $\chi$ such that $\chi = 1$ on an arbitrary compact set $K \subset \Omega_0$, we get

$$\|G(-\Delta_N)\Psi_\varepsilon\|^2_{L^2((0,T)\times K)} \leq \omega(\varepsilon) \left( \|\nabla \Psi_{0,\varepsilon}\|_2^2 + \|\varrho^{(1)}_{0,\varepsilon}\|_2^2 + \int_0^T \sum_{i=1}^5 \|G_{\varepsilon,i}\|_2^2 \right),$$

and finally

$$\|G(-\Delta_N)\Psi_\varepsilon\|^2_{L^2((0,T)\times K)} \to 0 \text{ as } \varepsilon \to 0. \tag{6.45}$$

## 7. Proof of the strong convergence

We are now able to conclude that (4.6) holds. Indeed, we have

$$\int_{\Omega_0} \widetilde{\mathbf{W}}_\varepsilon \cdot \varphi \, \mathrm{d}y = \int_{\Omega_0} \widetilde{\mathbf{W}}_\varepsilon \cdot \mathbf{H}(\varphi) \, \mathrm{d}y + \int_{\Omega_0} \widetilde{\mathbf{W}}_\varepsilon \cdot \mathbf{H}^\perp(\varphi) \, \mathrm{d}y \tag{7.1}$$

We already showed in section 5.2 that

$$\left[ t \to \int_{\Omega_0} \widetilde{\mathbf{W}}_\varepsilon \cdot \mathbf{H}(\varphi) \, \mathrm{d}y \right] \longrightarrow \left[ t \to \int_{\Omega_0} \widetilde{\mathbf{W}} \cdot \mathbf{H}(\varphi) \, \mathrm{d}y \right] \text{ strongly in } L^2(0,T). \tag{7.2}$$

Finally,

$$\begin{aligned} &\int_{\Omega_0} \widetilde{\mathbf{W}}_\varepsilon \cdot \mathbf{H}^\perp(\varphi) \, \mathrm{d}y = \int_{\Omega_0} \Psi_\varepsilon \, \mathrm{div}_y \, \mathbf{H}^\perp(\varphi) \, \mathrm{d}y \\ &= \int_{\Omega_0} G(-\Delta_N)\Psi_\varepsilon \, \mathrm{div}_y \, \mathbf{H}^\perp(\varphi) \, \mathrm{d}y + \int_{\Omega_0} (I - G(-\Delta_N))\Psi_\varepsilon \, \mathrm{div}_y \, \mathbf{H}^\perp(\varphi) \, \mathrm{d}y, \end{aligned} \tag{7.3}$$

where

$$\int_{\Omega_0} G(-\Delta_N)\Psi_\varepsilon \, \mathrm{div}_y \, \mathbf{H}^\perp(\varphi) \, \mathrm{d}y \longrightarrow 0,$$

due to (6.45) and

$$\int_{\Omega_0} (I - G(-\Delta_N))\Psi_\varepsilon \, \mathrm{div}_y \, \mathbf{H}^\perp(\varphi) \, \mathrm{d}y \longrightarrow 0,$$

since $I - G(-\Delta_N) \longrightarrow 0$ as $G \nearrow 1$ on $(0,\infty)$.

From (7.1), (7.2) and (7.3) it follows that

$$\left[ t \to \int_{\Omega_0} \widetilde{\mathbf{W}}_\varepsilon \cdot \varphi \, \mathrm{d}y \right] \longrightarrow \left[ t \to \int_{\Omega_0} \widetilde{\mathbf{W}} \cdot \varphi \, \mathrm{d}y \right] \text{ in } L^2(0,T).$$

Using change of variables and estimates (3.1), (3.4) we may conclude

$$\left[ t \to \int_{\Omega_t} \mathbf{u}_\varepsilon \cdot \varphi \, \mathrm{d}x \right] \longrightarrow \left[ t \to \int_{\Omega_t} \mathbf{U} \cdot \varphi \, \mathrm{d}x \right] \text{ in } L^2(0,T).$$

This, together with (4.3), implies (4.6).

## References

[1] Bogovskii, M.E.: *Solution of some vector analysis problems connected with operators div and grad* (in Russian) Trudy Sem. S.L. Sobolev, **80** (1980):5–40

[2] Cycon, H., L.; Froese, R., G.; Kirschi, W.; Simon, G: Schrödinger Operators with Application to Quantum Mechanics and Global Geometry, Springer (1987)

[3] Danchin, R.: *Low Mach number limit for viscous compressible flows* M2AN Math. Model Numer. Anal., **39** (2005) :459–475

[4] Desjardins, B., Grenier, E.: *Low Mach number limit of viscous compressible flows in the whole space* R. Soc. Lond. Proc. Ser. A Math. Phys. Eng. Sci., **455** (1999):2271–2279

[5] Desjardins, B., Grenier, E., Lions, P.-L., Masmoudi, N.: *Incompressible limit for solutions of the isentropic Navier-Stokes equations with Dirichlet boundary conditions* J. Math. Pures Appl., **78** (1999):461–471

[6] Feireisl, E.: *Flows of viscous compressible fluids under strong stratification: incompressible limits for long-range potential forces* Mathematical Models and Methods in Applied Sciences **21** (2003): 7-27

[7] Feireisl, E., Karper, T., Kreml, O., Stebel, J.: *Stability with respect to domain of the low Mach number limit of compressible viscous fluids* Mathematical Models and Methods in Applied Sciences **23**(2013): 2465-2493

[8] Feireisl, E., Kreml, O., Nečasová, Š., Neustupa, J., Stebel, J.: *Weak solutions to the barotropic Navier-Stokes system with slip boundary conditions in time dependent domains.* J. Differential Equations **254** (2013), no. 1, 125-140

[9] Feireisl, E., Kreml, O., Nečasová, Š., Neustupa, J., Stebel, J.: *Incompressible limits of fluids excited by moving boundaries.* SIAM J. Math. Anal. **46** (2014), No. 2, 1456-1471

[10] Feireisl, E., Novotný, A.: *Singular Limits in Thermodynamics of Viscous Fluids* Birkhaüser, Berlin (2009)

[11] Galdi, G.P.: *An introduction to the mathematical theory of the Navier - Stokes equations, I.* Springer-Verlag, New York (1994)

[12] Lighthill, J.: *On sound generated aerodynamically I. General theory* Proc. of the Royal Society of London, **A 211** (1952):564–587

[13] Lighthill, J.: *On sound generated aerodynamically II. General theory* Proc. of the Royal Society of London **A 222** (1954) :1–32

[14] Lions, P.-L.: *Mathematical topics in fluid dynamics, Vol.2, Compressible models* Oxford Science Publication, Oxford (1998)

[15] Lions, P.-L., Masmoudi, N.: *Incompressible Limit for a Viscous Compressible Fluid* J. Math. Pures Appl. **77** (1998): 585-627

[16] Masmoudi, N: *Examples of singular limits in hydrodynamics* In Handbook of Differential Equations, III, C. Dafermos, E. Feireisl Eds., Elsevier, Amsterdam (2006)

[17] Reed, M., Simon, B.: *Methods of modern mathematical physics. III. Scattering theory.* Academic Press [Harcourt Brace Jovanovich, Publishers], New York-London, 1979.

[18] Schochet, S.: *The mathematical theory of low Mach number flows.* M2ANMath. Model Numer. anal., **39** (2005), 441–458

[19] Sýkora, P.: *Compressible fluid motion in time dependent domains* Diploma Thesis, Charles University (2012)